\documentclass[11pt]{amsart}

\usepackage{amsfonts,amsmath,amssymb,amsthm}
\allowdisplaybreaks[4]

\usepackage{tikz}
\newcommand{\circled}[1]{%
  \tikz[baseline=(char.base)]{%
    \node[shape=circle,draw,inner sep=1pt] (char) {#1};%
  }%
}
\usepackage{braket}
\usepackage{graphicx}

\usepackage[colorlinks,linkcolor=red,citecolor=green]{hyperref}

\newtheorem{theorem}{Theorem}
\newtheorem{lemma}[theorem]{Lemma}

\numberwithin{theorem}{section}
\numberwithin{equation}{section}

\begin{document}

\setlength{\baselineskip}{1.1\baselineskip}
\title[A Liouville theorem for the $2$-Hessian equation on the Heisenberg group]{A Liouville theorem for the $2$-Hessian equation on the Heisenberg group}
\author{Wei Zhang}
\address{School of Mathematics and Statistics\\
         Lanzhou University\\
         Lanzhou, 730000, Gansu Province, China.}
\email{zhangw@lzu.edu.cn}
\author{Qi Zhou}
\address{School of Mathematics and Statistics\\
         Lanzhou University\\
         Lanzhou, 730000, Gansu Province, China.}
\email{zhouqimath20@lzu.edu.cn, zhouqi2025@lzu.edu.cn}
\maketitle

\begin{abstract}
 In this paper, we prove a Liouville theorem for the $2$-Hessian equation on the Heisenberg group $\mathbb{H}^n$. The result is obtained by choosing a suitable test function and using integration by parts to derive the necessary integral estimates.
\end{abstract}



\section{Introduction}

In 1981, the celebrated work of Gidas and Spruck \cite{Gidas-Spruck-1981} employed the integral estimate method introduced by Obata \cite{Obata-1971} to study the Liouville type properties of the equation
\begin{equation}\label{Equation-Laplace}
 \Delta u+u^\alpha=0\quad\mathrm{on}\ \mathbb{R}^n.
\end{equation}
They showed that the equation \eqref{Equation-Laplace} admits no positive entire solutions when $1<\alpha<\frac{n+2}{n-2}$. When $\alpha=\frac{n+2}{n-2}$, known as the critical exponent, Caffarelli, Gidas, and Spruck \cite{Caffarelli-Gidas-Spruck-1989} and Chen and Li \cite{Chen-Li-1991} classified all positive entire solutions by using the moving plane method. More recently, Ou \cite{Ou-2025} extended the integral estimate technique to obtain this classification in the lower dimensional cases $n=3,4,5$.

When the Laplacian operator in equation \eqref{Equation-Laplace} replaced by the $p$-Laplacian operator, resulting in the semilinear $p$-Laplace equation, the corresponding problem has also been extensively studied. Sciunzi \cite{Sciunzi-2016} and V\'{e}tois \cite{Vetois-2016} applied the moving plane method in their analysis. Integral estimates remain a powerful tool in this context, as demonstrated in the works of Serrin and Zou \cite{Serrin-Zou-2002}, Ciraolo, Figalli, and Roncoroni \cite{Ciraolo-Figalli-Roncoroni-2020}, Catino, Monticelli, and Roncoroni \cite{Catino-Monticelli-Roncoroni-2023}, and Ou \cite{Ou-2025}, among others.

For fully nonlinear elliptic operators, particularly the $k$-Hessian operator introduced by Caffarelli, Nirenberg, and Spruck \cite{Caffarelli-Nirenberg-Spruck-1985}, significant progress was initiated by Phuc and Verbitsky \cite{Phuc-Verbitsky-2006,Phuc-Verbitsky-2008}, who studied the $k$-Hessian equation
\begin{equation}\label{Equation-k-Hessian}
 \sigma_k(D^2u)=(-u)^\alpha\quad\mathrm{on}\ \mathbb{R}^n,
\end{equation}
and established nonexistence results for entire negative $k$-convex solutions under the conditions $2k<n$ and $k<\alpha\leq\frac{nk}{n-2k}$. They further showed that the exponent $\alpha=\frac{nk}{n-2k}$ is sharp for the inequality case. Here, $\sigma_k(D^2u)$ denotes the sum of all $k\times k$ principal minors of the Hessian matrix $D^2u$. A function $u\in C^2(\mathbb{R}^n)$ is called $k$-convex if it belongs to the cone
\[\Gamma_k=\{u\in C^2(\mathbb{R}^n) \mid \sigma_l(D^2u)\geq 0,\ l=1,2,\cdots,k\}.\]
The method used by Phuc and Verbitsky relies on the potential theory developed by Trudinger and Wang \cite{Trudinger-Wang-1997,Trudinger-Wang-1999,Trudinger-Wang-2002}. Subsequently, Ou \cite{Ou-2010} extended this result using integral estimate techniques when $\alpha\leq\frac{nk}{n-2k}$. Recently, Gao, Shi, and Zhang \cite{Gao-Shi-Zhang-2025} generalized these results to the $k$-$p$-Hessian equation introduced by Trudinger and Wang \cite{Trudinger-Wang-1999}. In the conformal setting, integral estimates continue to play an essential role in the study of Liouville type properties for $k$-Hessian equations, as seen in the works of Chang, Gursky, and Yang \cite{Chang-Gursky-Yang-2003}, Ou \cite{Ou-2013}, and so on.

We now turn our attention to the Heisenberg group $\mathbb{H}^n$. The equation corresponding to \eqref{Equation-Laplace} is given by
\begin{equation}\label{Equation-Laplace-Heisenberg}
 \Delta_{\mathbb{H}}u+u^\alpha=0\quad\mathrm{on}\ \mathbb{H}^n,
\end{equation}
where $\Delta_{\mathbb{H}}=\sum_{i=1}^{2n}X_iX_i$ denotes the sub-Laplacian on the Heisenberg group $\mathbb{H}^n$, and $\{X_i\}_{i=1,2,\cdots,2n}$ are the standard left-invariant vector fields. The details are provided in Section \ref{Section 2}. For $1<\alpha\leq\frac{Q}{Q-2}$, the nonexistence of positive entire solutions to \eqref{Equation-Laplace-Heisenberg} was established by Birindelli, Dolcetta, and Cutri \cite{Birindelli-Capuzzo-Cutri-1997} using integral estimates. They also showed that the exponent $\frac{Q}{Q-2}$ is sharp for the inequality version. Later, Xu \cite{Xu-2009} extended this result to the broader range $1<\alpha<\frac{Q(Q+2)}{(Q-1)^2}$. In a notable contribution \cite{Ma-Ou-2023}, Ma and Ou developed a key identity originally due to Jerison and Lee \cite{Jerison-Lee-1988}, and applied integration by parts techniques to obtain a complete classification of all positive solutions when $1<\alpha<\frac{Q+2}{Q-2}$. In the parabolic case, the corresponding result was proved by Wei and Wu \cite{Wei-Wu-2024} for $1<\alpha<1+\frac{4Q-10}{Q^2-4Q+6}$.

The critical case $\alpha=\frac{Q+2}{Q-2}$ is closely related to the CR Yamabe problem and the sharp Sobolev inequality on $\mathbb{H}^n$, as studied extensively by Jerison and Lee in their series of works \cite{Jerison-Lee-1987,Jerison-Lee-1988,Jerison-Lee-1989}. In a recent development, Catino, Li, Monticelli, and Roncoroni \cite{Catino-Li-Monticelli-Roncoroni-2023} adopted an approach inspired by Ou \cite{Ou-2025}, who studied the $p$-Laplace equation on $\mathbb{R}^n$, to give a complete classification of entire positive solutions of \eqref{Equation-Laplace-Heisenberg} in the case of $\mathbb{H}^1$ without additional assumptions. In higher dimensions, they and Flynn and V\'{e}tois \cite{Flynn-Vetois-2023} each treated the problem under different assumptions.

It is natural to ask whether equation \eqref{Equation-k-Hessian}, in the Heisenberg group setting, admits only trivial $k$-convex nonpositive solutions. This problem corresponds to extending the nonexistence results of Phuc and Verbitsky \cite{Phuc-Verbitsky-2006,Phuc-Verbitsky-2008} and Ou \cite{Ou-2010} from Euclidean space to the Heisenberg group, and can also be interpreted as a generalization of the works of Birindelli, Dolcetta, and Cutri \cite{Birindelli-Capuzzo-Cutri-1997}, and Ma and Ou \cite{Ma-Ou-2023} to the $k$-Hessian framework.

Before presenting our main result, we briefly recall the notion of convexity for functions on the Heisenberg group. Various definitions of convexity in this setting have been proposed and compared, including horizontal convexity introduced by Danielli, Garofalo, and Nhieu \cite{Danielli-Garofalo-Nhieu-2003}, and viscosity convexity by Lu, Manfredi, and Stroffolini \cite{Lu-Manfredi-Stroffolini-2004}. These different notions have since been shown to be equivalent. In this paper, a function $u\in C^2(\mathbb{H}^n)$ is said to be convex if the symmetric matrix
\[\mathrm{Hess}_Xu=\left[\frac{X_iX_ju+X_jX_iu}{2}\right]_{i,j=1,2,\cdots,2n}\]
is semi-positive definite in $\mathbb{H}^n$. As in the Euclidean setting, we say that $u\in C^2(\mathbb{H}^n)$ is $k$-convex if $u\in\Gamma_k$, where the cone $\Gamma_k$ is defined by
\[\Gamma_k=\{u\in C^2(\mathbb{H}^n)\mid\sigma_l(\mathrm{Hess}_Xu)\geq 0,\ l=1,2,\cdots,k\}.\]
This paper is devoted to the case $k=2$, in which we establish a Liouville type theorem for $2$-convex solutions of the following equation
\begin{equation}\label{Equation-sigma_2}
 \sigma_2(\mathrm{Hess}_Xu)=(-u)^\alpha\quad\mathrm{on}\ \mathbb{H}^n.
\end{equation}

We now state our main result.
\begin{theorem}\label{Main Theorem}
 If $\alpha\leq\frac{2Q}{Q-4}$ with $Q>4$, then \eqref{Equation-sigma_2} has no negative solution in $\Gamma_2$.
\end{theorem}

Our proof mainly relies on integral estimate techniques. It is worth noting that, since the Heisenberg group $\mathbb{H}^n$ is a noncommutative Lie group, integration by parts inevitably gives rise to additional commutator terms, which make the computations substantially more involved than in the Euclidean setting.

This paper is organized as follows. In Section \ref{Section 2}, we review basic definitions of the Heisenberg group and collect several results from linear algebra that will be frequently used later. Section \ref{Section 3} is devoted to establishing Theorem \ref{Main Theorem}.

\section{Preliminaries}\label{Section 2}
The Heisenberg group $\mathbb{H}^n$ can be represented as the Euclidean space $\mathbb{R}^n\times\mathbb{R}^n\times\mathbb{R}$, endowed with a specific group multiplication law. For two elements $g=(x,y,t)$ and $g_0=(x_0,y_0,t_0)$ in $\mathbb{H}^n$, the multiplication is defined as
\[g\circ g_0=\big(x+x_0,y+y_0,t+t_0+2(\Braket{y,x_0}-\Braket{x,y_0})\big),\]
where $\Braket{\cdot,\cdot}$ denotes the standard inner product in $\mathbb{R}^n$. For indices $i,j=1,2,\cdots,n$, the left-invariant vector fields $X_i$ and $X_{n+j}$ in $\mathbb{H}^n$ are constructed such that they coincide with $\frac{\partial}{\partial x_i}$ and $\frac{\partial}{\partial y_j}$, respectively, at the origin. Explicitly, these vector fields are given by
\[X_i=\frac{\partial}{\partial x_i}+2y_i\frac{\partial}{\partial t},\quad X_{n+j}=\frac{\partial}{\partial y_j}-2x_j\frac{\partial}{\partial t}.\]
Additionally, the left-invariant vector field $T$, which corresponds to to $\frac{\partial}{\partial t}$ at the origin, is defined as
\[T=\frac{\partial}{\partial t}.\]
Together, the $2n+1$ vector fields $X_1, X_2,\cdots,X_{2n},T$ form a basis for the Lie algebra of $\mathbb{H}^n$. The commutation relations among these vector fields are given by
\begin{equation}\label{Relation-1}
 [X_i,X_{n+j}]=-4\delta_{ij}T,\quad \forall\ i,j=1,2,\cdots,n,
\end{equation}
where $\delta_{ij}$ is the Kronecker symbol. All other commutators between these vector fields vanish, i.e.,
\begin{equation}\label{Relation-2}
 [X_i,X_j]=0,\quad [X_{n+i},X_{n+j}]=0,\quad\forall\ i,j=1,2,\cdots,n.
\end{equation}
This structure defines the noncommutative nature of the Heisenberg group.

Let us define the gauge function
\[\rho(g)=\big((|x|^2+|y|^2)^2+t^2\big)^{1/4},\]
where $g=(x,y,t)$ represents a point in $\mathbb{H}^n$. This gauge function allows us to introduce a distance function $\mathrm{d}(g,g_0)$ between two points $g$ and $g_0$ in $\mathbb{H}^n$, given by
\[\mathrm{d}(g,g_0)=\rho(g_0^{-1}\circ g).\]
The Heisenberg group $\mathbb{H}^n$ also admits dilations $\delta_\lambda$ for $\lambda>0$, defined as
\[\delta_\lambda(g)=(\lambda x,\lambda y,\lambda^2t).\] 
The distance function $\mathrm{d}$ is homogeneous of degree one under these dilations, meaning
\[\mathrm{d}(\delta_\lambda g,\delta_\lambda g_0)=\lambda\mathrm{d}(g,g_0)\]
for any $\lambda>0$ and $g, g_0\in\mathbb{H}^n$. Using this distance, we can define a ball $B_R(g_0)$ centered at $g_0$ with radius $R$ as
\[B_R(g_0)=\{g\in\mathbb{H}^n \mid \mathrm{d}(g,g_0)<R\}.\]
The Lebesgue measure of such a ball scales according to
\[|B_R(g_0)|=R^Q|B_1(0)|,\]
where $Q=2n+2$ is the homogeneous dimension of $\mathbb{H}^n$. For further details on the Heisenberg group and its properties, Stein's book \cite{Stein-1993} provides a comprehensive reference.

As a preparation for the proof of our main theorem, we now collect some linear algebra preliminaries.

Let $A=(a_{ij})_{n\times n}$ be a symmetric matrix. For notational convenience, we denote the first and second derivatives of the second  elementary symmetric function $\sigma_2(A)$ with respect to the entries of $A$ by
\[\sigma_2^{ij}=\frac{\partial\sigma_2(A)}{\partial a_{ij}}\quad\mathrm{and}\ \sigma_2^{ij,kl}=\frac{\partial^2\sigma_2(A)}{\partial a_{ij}\partial a_{kl}}.\]
It follows directly from the homogeneity of $\sigma_2$ that \[\sigma_2(A)=\frac{1}{2}\sum_{i,j=1}^n\sigma_2^{ij}a_{ij}.\] 
Since $\sigma_2(A)$ denotes the sum of all $2\times 2$ principal minors of $A$, it can be explicitly written as
\[\sigma_2(A)=\frac{1}{2}\sum_{i\neq j}(a_{ii}a_{jj}-a_{ij}a_{ji}).\]
For this expression, we readily deduce that the second derivatives satisfy $\sigma_2^{ii,jj}=1$ and $\sigma_2^{ij,ji}=-1$ for $i\neq j$, while all other second order partial derivatives vanish. With this observation, we can easily derive the following lemmas, originally proved by Trudinger and the first author \cite{Trudinger-Zhang-2013}.

\begin{lemma}[Lemma 2.1 in \cite{Trudinger-Zhang-2013}]\label{Lemma 2.1}
 Let $A=(a_{ij})_{n\times n}$ be a symmetric matrix, then
 \[\frac{\partial^2\sigma_2(A)}{\partial a_{ij}\partial a_{kl}}=-\frac{\partial^2\sigma_2(A)}{\partial a_{kj}\partial a_{il}}.\]
\end{lemma}
\begin{lemma}[Lemma 2.2 in \cite{Trudinger-Zhang-2013}]\label{Lemma 2.2}
 Let $A=(a_{ij})_{n\times n}$ be a symmetric matrix, then we have the following commutation formula
 \[\frac{\partial^2\sigma_2(A)}{\partial a_{ij}\partial a_{kl}}-\frac{\partial^2\sigma_2(A)}{\partial a_{ij}\partial a_{lk}}=\frac{\partial^2\sigma_2(A)}{\partial a_{ik}\partial a_{jl}}.\]
\end{lemma}

\section{Proof of Theorem \ref{Main Theorem}}\label{Section 3}

In this section, we prove Theorem \ref{Main Theorem} using integral estimates. Our strategy involves multiplying both sides of equation \eqref{Equation-sigma_2} by $(-u)^\delta\eta^\theta$, followed by integration by parts and repeated applications of Young's inequality to derive the desired estimate. The $2$-Hessian operator is a fully nonlinear elliptic operator, which contributes to the analytical complexity of the equation. Moreover, since the Heisenberg group is a noncommutative Lie group, integration by parts leads to the appearance of commutator terms when the order of derivatives is exchanged. This feature increases the difficulty to the calculation.

\begin{proof}[Proof of Theorem \ref{Main Theorem}]
 Assume $u<0$ be a solution of \eqref{Equation-sigma_2} in $\Gamma_2$. Let $\eta$ be a  real smooth cut off function satisfying
 \[\begin{cases}
  \eta\equiv 1 & \mathrm{in}\ B_R,\\
  0\leq\eta\leq 1 & \mathrm{in}\ B_{2R},\\
  \eta\equiv 0 & \mathrm{in}\ \mathbb{H}^n\setminus B_{2R},\\
  |\partial\eta|\leq\frac{C(Q)}{R} & \mathrm{on}\ \mathbb{H}^n,
 \end{cases}\]
 where $|\partial\eta|^2=\sum_{i=1}^{2n}(X_i\eta)^2$. Here, $C(Q)$ means a positive constant depending only on the homogeneous dimension $Q$ of the Heisenberg group $\mathbb{H}^n$.
 
 Multiplying $2(-u)^\delta\eta^\theta$ on both side of equation \eqref{Equation-sigma_2}, and integrating on $\mathbb{H}^n$ (we always omit the volume form), we have
 \begin{equation}\label{Inequality-1}
  2\int(-u)^{\alpha+\delta}\eta^\theta= 2\int\sigma_2(\mathrm{Hess}_Xu)(-u)^\delta\eta^\theta,
 \end{equation}
 where $\delta$ and $\theta$ are constants to be determined later. Note that
 \[\sigma_2(\mathrm{Hess}_Xu)=\frac{1}{2}\sum_{i,j=1}^{2n}\sigma_2^{ij}\frac{X_iX_ju+X_jX_iu}{2}.\]
 In the following, we write $\sigma_2(\mathrm{Hess}_Xu)$ simply as $\sigma_2$. For the right hand side of \eqref{Inequality-1}, integrating by parts, we obtain
 \begin{equation}\label{Equality-1+2+3-(1)}
  \begin{aligned}
   2\int\sigma_2(-u)^\delta\eta^\theta
   =&\sum_{i,j=1}^{2n}\int\sigma_2^{ij}\frac{X_iX_ju+X_jX_iu}{2}(-u)^\delta\eta^\theta\\
   =&\sum_{i,j=1}^{2n}\int\sigma_2^{ij}X_iX_ju(-u)^\delta\eta^\theta\\
   =&-\sum_{i,j=1}^{2n}\int X_i(\sigma_2^{ij})X_ju(-u)^\delta\eta^\theta+\delta\sum_{i,j=1}^{2n}\int\sigma_2^{ij}X_iuX_ju(-u)^{\delta-1}\eta^\theta\\
   &-\theta\sum_{i,j=1}^{2n}\int\sigma_2^{ij}X_juX_i\eta(-u)^\delta\eta^{\theta-1}\\
   \stackrel{\triangle}{=}&~\circled{1}+\circled{2}+\circled{3}.
  \end{aligned}
 \end{equation}
 
 We perform the following calculations using the Einstein summation convention, in which repeated indices denote summation. Unless otherwise specified, the indices $i,j,k,l$ range from $1$ to $2n$.
 
 Let us start by computing terms \circled{1}, \circled{2}, and \circled{3} in detail. For the term \circled{1}, we permute the indices and apply Lemma \ref{Lemma 2.1} to obtain
 \begin{align*}
  \circled{1}
  =&-\int\sigma_2^{ij,kl}X_i\left(\frac{X_kX_lu+X_lX_ku}{2}\right)X_ju(-u)^\delta\eta^\theta\\
  =&-\frac{1}{2}\int\sigma_2^{ij,kl}X_iX_kX_luX_ju(-u)^\delta\eta^\theta-\frac{1}{2}\int\sigma_2^{ij,kl}X_iX_lX_kuX_ju(-u)^\delta\eta^\theta\\
  =&-\frac{1}{2}\int\sigma_2^{ij,kl}(X_iX_kX_lu-X_kX_lX_iu)X_ju(-u)^\delta\eta^\theta.
 \end{align*}
 Since
 \begin{align*}
  X_iX_kX_lu-X_kX_lX_iu
  =&~[X_i,X_k]X_lu+X_kX_iX_lu-X_kX_lX_iu\\
  =&~[X_i,X_k]X_lu+X_k[X_i,X_l]u,
 \end{align*}
 we can rewrite \circled{1} as
 \begin{align*}
  \circled{1}
  =&-\frac{1}{2}\int\sigma_2^{ij,kl}[X_i,X_k]X_luX_ju(-u)^\delta\eta^\theta-\frac{1}{2}\int\sigma_2^{ij,kl}X_k[X_i,X_l]uX_ju(-u)^\delta\eta^\theta\\
  \stackrel{\triangle}{=}&~\circled{1}_1+\circled{1}_2.
 \end{align*}
 By interchanging the indices $i$ and $k$ in the term $\circled{1}_2$, and using Lemmas \ref{Lemma 2.1} and \ref{Lemma 2.2}, we deduce
 \begin{align*}
  \circled{1}_2
  =&~\frac{1}{2}\int\sigma_2^{ij,kl}X_i[X_k,X_l]uX_ju(-u)^\delta\eta^\theta\\
  =&~\frac{1}{4}\int(\sigma_2^{ij,kl}-\sigma_2^{ij,lk})X_i[X_k,X_l]uX_ju(-u)^\delta\eta^\theta\\
  =&~\frac{1}{4}\int\sigma_2^{ik,jl}X_i[X_k,X_l]uX_ju(-u)^\delta\eta^\theta\\
  =&~\frac{1}{4}\int\sigma_2^{lk,ji}X_l[X_k,X_i]uX_ju(-u)^\delta\eta^\theta\\
  =&~-\frac{1}{4}\int\sigma_2^{ij,kl}[X_i,X_k]X_luX_ju(-u)^\delta\eta^\theta\\
  =&~\frac{1}{2}~\circled{1}_1.
 \end{align*}
 We thus obtain
 \begin{align*}
  \circled{1}=-\frac{3}{4}\int\sigma_2^{ij,kl}[X_i,X_k]X_luX_ju(-u)^\delta\eta^\theta.
 \end{align*}
 To further evaluate term \circled{1}, we carry out integration by parts, which yields
 \begin{align*}
  \circled{1}
  =&~\frac{3}{4}\int X_l(\sigma_2^{ij,kl})[X_i,X_k]uX_ju(-u)^\delta\eta^\theta+\frac{3}{4}\int\sigma_2^{ij,kl}[X_i,X_k]uX_lX_ju(-u)^\delta\eta^\theta\\
  &-\frac{3}{4}\delta\int\sigma_2^{ij,kl}[X_i,X_k]uX_juX_lu(-u)^{\delta-1}\eta^\theta+\frac{3}{4}\theta\int\sigma_2^{ij,kl}[X_i,X_k]uX_juX_l\eta(-u)^\delta\eta^{\theta-1}.
 \end{align*}
 The third term in the right hand side of the above equality is vanishing by the skew symmetric for $j$ and $l$. It follows that
 \begin{align*}
  \circled{1}
  =&-\frac{3}{8}\int\sigma_2^{ij,kl}[X_i,X_k]u[X_j,X_l]u(-u)^\delta\eta^\theta+\frac{3}{4}\theta\int\sigma_2^{ij,kl}[X_i,X_k]uX_juX_l\eta(-u)^\delta\eta^{\theta-1}\\
  \stackrel{\triangle}{=}&~\circled{1}_3+\circled{1}_4.
 \end{align*}
 Using the commutation relations \eqref{Relation-1} and \eqref{Relation-2}, from the facts $\sigma_2^{ii,jj}=1$ and $\sigma_2^{ij,ji}=-1$ for $i\neq j$, with all other terms vanishing,
 \begin{align*}
  \circled{1}_3
  =&-\frac{3}{8}\cdot 4\sum_{i<k,j<l}\int\sigma_2^{ij,kl}[X_i,X_k]u[X_j,X_l]u(-u)^\delta\eta^\theta\\
  =&-\frac{3}{2}\sum_{i,j=1}^n\int\sigma_2^{ij,n+in+j}\cdot 16(Tu)^2(-u)^\delta\eta^\theta\\
  =&-24\sum_{i=1}^n\int\sigma_2^{ii,n+in+i}(Tu)^2(-u)^\delta\eta^\theta\\
  =&-24n\int(Tu)^2(-u)^\delta\eta^\theta
 \end{align*}
 and the second term can be derived as
 \begin{align*}
  \circled{1}_4
  =&~\frac{3}{4}\theta\cdot 2\sum_{i<k}\sum_{j,l=1}^{2n}\int\sigma_2^{ij,kl}[X_i,X_k]uX_juX_l\eta(-u)^\delta\eta^{\theta-1}\\
  =&-\frac{3}{2}\theta\sum_{i=1}^n\sum_{j,l=1}^{2n}\int\sigma_2^{ij,n+il}\cdot 4(Tu)X_juX_l\eta(-u)^\delta \eta^{\theta-1}\\
  =&~12\theta\sum_{i=1}^n\int(Tu)X_{n+i}uX_i\eta(-u)^\delta\eta^{\theta-1},
 \end{align*}
 which together imply
 \[\circled{1}=-24n\int(Tu)^2(-u)^\delta\eta^\theta+12\theta\sum_{i=1}^n\int(Tu)X_{n+i}uX_i\eta(-u)^\delta\eta^{\theta-1}.\]
 With the help of the Cauchy-Schwarz inequality and the known properties of the cut-off function $\eta$, we can estimate
 \begin{equation}\label{Estimate-1}
  \circled{1}\leq-24n\int(Tu)^2(-u)^\delta\eta^\theta+12|\theta|\frac{C(Q)}{R}\int|Tu||\partial u|(-u)^\delta\eta^{\theta-1}.
 \end{equation}
 
 For the term \circled{2}, we have
 \begin{align*}
  \circled{2}
  =&~\delta\int\sigma_2^{ij}X_iuX_ju(-u)^{\delta-1}\eta^\theta\\
  =&~\delta\int\left(\sigma_1\delta_{ij}-\frac{X_jX_iu+X_iX_ju}{2}\right)X_iuX_ju(-u)^{\delta-1}\eta^\theta\\
  =&~\delta\int\sigma_1|\partial u|^2(-u)^{\delta-1}\eta^\theta-\delta\int X_iX_juX_iuX_ju(-u)^{\delta-1}\eta^\theta,
 \end{align*}
 where $|\partial u|^2=\sum_{i=1}^{2n}(X_iu)^2$. Let
 \[\circled{2}_1=-\delta\int X_iX_juX_iuX_ju(-u)^{\delta-1}\eta^\theta.\]
 By applying integration by parts, we derive
 \begin{align*}
  \circled{2}_1
  =&~\delta\int \sigma_1|\partial u|^2(-u)^{\delta-1}\eta^\theta-\circled{2}_1-\delta(\delta-1)\int|\partial u|^4(-u)^{\delta-2}\eta^\theta\\
  &+\delta\theta\int|\partial u|^2X_iuX_i\eta(-u)^{\delta-1}\eta^{\theta-1},
 \end{align*}
 which leads to
 \begin{align*}
  \circled{2}_1
  =&~\frac{1}{2}\delta\int\sigma_1|\partial u|^2(-u)^{\delta-1}\eta^\theta-\frac{1}{2}\delta(\delta-1)\int|\partial u|^4(-u)^{\delta-2}\eta^\theta\\
  &+\frac{1}{2}\delta\theta\int|\partial u|^2X_iuX_i\eta(-u)^{\delta-1}\eta^{\theta-1}.
 \end{align*}
 As a result, it can be concluded that
 \begin{align*}
  \circled{2}
  =&~\frac{3}{2}\delta\int\sigma_1|\partial u|^2(-u)^{\delta-1}\eta^\theta-\frac{1}{2}\delta(\delta-1)\int|\partial u|^4(-u)^{\delta-2}\eta^\theta\\
  &+\frac{1}{2}\delta\theta\int|\partial u|^2X_iuX_i\eta(-u)^{\delta-1}\eta^{\theta-1},
 \end{align*}
 and by the Cauchy-Schwarz inequality, this gives
 \begin{equation}\label{Estimate-2}
  \begin{aligned}
   \circled{2}
   \leq&~\frac{3}{2}\delta\int\sigma_1|\partial u|^2(-u)^{\delta-1}\eta^\theta-\frac{1}{2}\delta(\delta-1)\int|\partial u|^4(-u)^{\delta-2}\eta^\theta\\
   &+\frac{1}{2}|\delta||\theta|\frac{C(Q)}{R}\int|\partial u|^3(-u)^{\delta-1}\eta^{\theta-1}.
  \end{aligned}
 \end{equation}
 
 We now turn to term \circled{3}. Following the same reasoning as for term \circled{2}, it is straightforward to verify that
 \begin{align*}
  \circled{3}
  =&-\theta\int\left(\sigma_1\delta_{ij}-\frac{X_jX_iu+X_iX_ju}{2}\right)X_juX_i\eta(-u)^\delta\eta^{\theta-1}\\
  =&-\theta\int\sigma_1X_iuX_i\eta(-u)^\delta\eta^{\theta-1}\\
  &+\frac{1}{2}\theta\int X_jX_iuX_juX_i\eta(-u)^\delta\eta^{\theta-1}+\frac{1}{2}\theta\int X_iX_juX_juX_i\eta(-u)^\delta\eta^{\theta-1}.
 \end{align*}
 We label the last two terms by $\circled{3}_1$ and $\circled{3}_2$, respectively. Using integration by parts to both, we find
 \begin{align*}
  \circled{3}_1
  =&-\frac{1}{2}\theta\int\sigma_1X_iuX_i\eta(-u)^\delta\eta^{\theta-1}-\frac{1}{2}\theta\int X_iuX_juX_jX_i\eta(-u)^\delta\eta^{\theta-1}\\
  &+\frac{1}{2}\delta\theta\int|\partial u|^2X_iuX_i\eta(-u)^{\delta-1}\eta^{\theta-1}-\frac{1}{2}\theta(\theta-1)\int X_iuX_juX_i\eta X_j\eta(-u)^\delta\eta^{\theta-2},
 \end{align*}
 \begin{align*}
  \circled{3}_2
  =&-\circled{3}_2-\frac{1}{2}\theta\int|\partial u|^2\sigma_1(\mathrm{Hess}_X\eta)(-u)^\delta\eta^{\theta-1}\\
  &+\frac{1}{2}\delta\theta\int|\partial u|^2X_iuX_i\eta(-u)^{\delta-1}\eta^{\theta-1}-\frac{1}{2}\theta(\theta-1)\int|\partial u|^2|\partial\eta|^2(-u)^\delta\eta^{\theta-2},
 \end{align*}
 so that
 \begin{align*}
  \circled{3}_2
  =&-\frac{1}{4}\theta\int|\partial u|^2\sigma_1(\mathrm{Hess}_X\eta)(-u)^\delta\eta^{\theta-1}+\frac{1}{4}\delta\theta\int|\partial u|^2X_iuX_i\eta(-u)^{\delta-1}\eta^{\theta-1}\\
  &-\frac{1}{4}\theta(\theta-1)\int|\partial u|^2|\partial\eta|^2(-u)^\delta\eta^{\theta-2},
 \end{align*}
 and therefore
 \begin{align*}
  \circled{3}
  =&-\frac{3}{2}\theta\int\sigma_1X_iuX_i\eta(-u)^\delta\eta^{\theta-1}+\frac{3}{4}\delta\theta\int|\partial u|^2X_iuX_i\eta(-u)^{\delta-1}\eta^{\theta-1}\\
  &-\frac{1}{2}\theta(\theta-1)\int X_iuX_juX_i\eta X_j\eta(-u)^\delta\eta^{\theta-2}-\frac{1}{2}\theta\int X_iuX_juX_jX_i\eta(-u)^\delta\eta^{\theta-1}\\
  &-\frac{1}{4}\theta(\theta-1)\int|\partial u|^2|\partial\eta|^2(-u)^\delta\eta^{\theta-2}-\frac{1}{4}\theta\int|\partial u|^2\sigma_1(\mathrm{Hess}_X\eta)(-u)^\delta\eta^{\theta-1}.
 \end{align*}
 Since $\eta\leq 1$ throughout $\mathbb{H}^n$, this yields
 \begin{equation}\label{Estimate-3}
  \begin{aligned}
   \circled{3}
   \leq&~\frac{3}{2}|\theta|\frac{C(Q)}{R}\int\sigma_1|\partial u|(-u)^\delta\eta^{\theta-1}+\frac{3}{4}|\delta||\theta|\frac{C(Q)}{R}\int|\partial u|^3(-u)^{\delta-1}\eta^{\theta-1}\\
   &+\frac{3}{4}|\theta|(|\theta-1|+1)\frac{C(Q)}{R^2}\int|\partial u|^2(-u)^\delta\eta^{\theta-2}.
  \end{aligned}
 \end{equation}
 From \eqref{Equality-1+2+3-(1)}, together with \eqref{Estimate-1}, \eqref{Estimate-2}, and \eqref{Estimate-3}, we derive the integral estimate
 \begin{equation}\label{Inequality-1+2+3}
  \begin{aligned}
   2\int\sigma_2(-u)^\delta\eta^\theta
   \leq&-24n\int(Tu)^2(-u)^\delta\eta^\theta+\frac{3}{2}\delta\int\sigma_1|\partial u|^2(-u)^{\delta-1}\eta^\theta\\
   &-\frac{1}{2}\delta(\delta-1)\int|\partial u|^4(-u)^{\delta-2}\eta^\theta\\
   &+\frac{3}{4}|\theta|(|\theta-1|+1)\frac{C(Q)}{R^2}\int|\partial u|^2(-u)^\delta\eta^{\theta-2}\\
   &+12|\theta|\frac{C(Q)}{R}\int|Tu||\partial u|(-u)^\delta\eta^{\theta-1}\\
   &+\frac{3}{2}|\theta|\frac{C(Q)}{R}\int\sigma_1|\partial u|(-u)^\delta\eta^{\theta-1}\\
   &+\frac{5}{4}|\delta||\theta|\frac{C(Q)}{R}\int|\partial u|^3(-u)^{\delta-1}\eta^{\theta-1}.
  \end{aligned}
 \end{equation}
 
 In what follows, we estimate the last four terms in \eqref{Inequality-1+2+3}, labeled \circled{4}, \circled{5}, \circled{6}, and \circled{7}. The key idea is to use Young's inequality to express them using other existing terms in \eqref{Inequality-1+2+3}. As the coefficients involved are nonnegative, we will omit them in the applications of Young's inequality that follows, without affecting the validity of the estimates.
 
 For the term \circled{4}, by applying Young's inequality with the exponent pair $(2,2)$, we have
 \begin{equation}\label{Estimate-4}
  \circled{4}\leq\varepsilon_1\int|\partial u|^4(-u)^{\delta-2}\eta^\theta+C(\varepsilon_1)R^{-4}\int(-u)^{\delta+2}\eta^{\theta-4}.
 \end{equation}
 In the same manner, the term \circled{5} can be estimated as
 \begin{equation}\label{Estimate-5}
  \begin{aligned}
   \circled{5}
   \leq&~\varepsilon_2\int(Tu)^2(-u)^\delta\eta^\theta+C(\varepsilon_2)R^{-2}\int|\partial u|^2(-u)^\delta\eta^{\theta-2}\\
   \leq&~\varepsilon_2\int(Tu)^2(-u)^\delta\eta^\theta+\varepsilon_1\int|\partial u|^4(-u)^{\delta-2}\eta^\theta\\
   &+C(\varepsilon_1,\varepsilon_2)R^{-4}\int(-u)^{\delta+2}\eta^{\theta-4}.
  \end{aligned}
 \end{equation}
 
 For the term \circled{6}, it is clear that
 \[\circled{6}\leq\varepsilon_3\int\sigma_1|\partial u|^2(-u)^{\delta-1}\eta^\theta+C(\varepsilon_3)R^{-2}\int\sigma_1(-u)^{\delta+1}\eta^{\theta-2}.\]
 Denote the last term in this inequality by $\circled{6}_1$. By carrying out integration by parts and successively applying Young's inequality with the exponent pairs $(2,2)$ and $(4,\frac{4}{3})$, this procedure yields
 \begin{align*}
  \circled{6}_1
  =&~(\delta+1)C(\varepsilon_3)R^{-2}\int|\partial u|^2(-u)^\delta\eta^{\theta-2}-(\theta-2)C(\varepsilon_3)R^{-2}\int X_iuX_i\eta(-u)^{\delta+1}\eta^{\theta-3}\\
  \leq&~|\delta+1|C(\varepsilon_3)R^{-2}\int|\partial u|^2(-u)^\delta\eta^{\theta-2}+|2-\theta|C(\varepsilon_3)R^{-3}\int|\partial u|(-u)^{\delta+1}\eta^{\theta-3}\\
  \leq&~(\varepsilon_1+\varepsilon_4)\int|\partial u|^4(-u)^{\delta-2}\eta^\theta+C(\varepsilon_1,\varepsilon_3,\varepsilon_4)R^{-4}\int(-u)^{\delta+2}\eta^{\theta-4}.
 \end{align*}
 Hence, an estimate for \circled{6} is given by
 \begin{equation}\label{Estimate-6}
  \begin{aligned}
   \circled{6}
   \leq&~\varepsilon_3\int\sigma_1|\partial u|^2(-u)^{\delta-1}\eta^\theta+(\varepsilon_1+\varepsilon_4)\int|\partial u|^4(-u)^{\delta-2}\eta^\theta\\
   &+C(\varepsilon_1,\varepsilon_3,\varepsilon_4)R^{-4}\int(-u)^{\delta+2}\eta^{\theta-4}.
  \end{aligned}
 \end{equation}
 
 We now turn on the final term \circled{7}, which can be estimated directly using Young’s inequality with the exponent pair $(\frac{4}{3},4)$, resulting in
 \begin{equation}\label{Estimate-7}
  \circled{7}\leq\varepsilon_5\int|\partial u|^4(-u)^{\delta-2}\eta^\theta+C(\varepsilon_5)R^{-4}\int(-u)^{\delta+2}\eta^{\theta-4}.
 \end{equation}
 
 By inserting the bounds from \eqref{Estimate-4}, \eqref{Estimate-5}, \eqref{Estimate-6}, and \eqref{Estimate-7} into \eqref{Inequality-1+2+3}, we arrive at the integral estimate
 \begin{equation}\label{Inequality-2}
 \begin{aligned}
  &~2\int\sigma_2(-u)^\delta\eta^\theta\\
  \leq&~(\varepsilon_2-24n)\int(Tu)^2(-u)^\delta\eta^\theta+\left(\varepsilon_3+\frac{3}{2}\delta\right)\int\sigma_1|\partial u|^2(-u)^{\delta-1}\eta^\theta\\
  &+\left(\varepsilon_1+\varepsilon_4+\varepsilon_5-\frac{1}{2}\delta(\delta-1)\right)\int|\partial u|^4(-u)^{\delta-2}\eta^\theta\\
  &+C(\varepsilon_1,\varepsilon_2,\varepsilon_3,\varepsilon_4,\varepsilon_5)R^{-4}\int(-u)^{\delta+2}\eta^{\theta-4}.
 \end{aligned}
 \end{equation}
 We choose $\varepsilon_i$ ($i=1,2,\cdots,5$) to be sufficiently small and fix them, and at the same time take $\delta<0$, so as to ensure that the coefficients of the first three terms on the right hand side of \eqref{Inequality-2} are nonpositive. Note that $u<0$ and in $\Gamma_2$. Consequently, based on \eqref{Inequality-1}, we have
 \begin{equation}\label{Inequality-3}
  2\int(-u)^{\alpha+\delta}\eta^\theta\leq C(Q)R^{-4}\int(-u)^{\delta+2}\eta^{\theta-4}.
 \end{equation}
 
 Observe that, by employing Young's inequality with the exponent pair $(\frac{\alpha+\delta}{2+\delta},\frac{\alpha+\delta}{\alpha-2})$ on \eqref{Inequality-3}, we get
 \[2\int(-u)^{\alpha+\delta}\eta^\theta\leq\varepsilon_6\int(-u)^{\alpha+\delta}\eta^\theta+C(\varepsilon_6)R^{-4\frac{\alpha+\delta}{\alpha-2}}\int\eta^{\theta-4\frac{\alpha+\delta}{\alpha-2}},\]
 from which it follows that
 \begin{equation}\label{Inequality-4}
  (2-\varepsilon_6)\int(-u)^{\alpha+\delta}\eta^\theta\leq C(\varepsilon_6)R^{Q-4\frac{\alpha+\delta}{\alpha-2}}.
 \end{equation}
 Here, it is necessary that $\frac{\alpha+\delta}{2+\delta}>1$, $\alpha\neq-\delta\neq 2$ and that $\theta>0$ is sufficiently large so that $\theta-4\frac{\alpha+\delta}{\alpha-2}>0$. If we require that $Q-4\frac{\alpha+\delta}{\alpha-2}<0$ and choose $\varepsilon_6$ small enough, then from \eqref{Inequality-4}, letting $R\rightarrow+\infty$, we obtain $u\equiv 0$ in $\mathbb{H}^n$, which contradicts our assumption that $u<0$. Now, based on the conditions $\frac{\alpha+\delta}{2+\delta}>1$, $\alpha\neq-\delta\neq 2$, $\delta<0$ and $Q-4\frac{\alpha+\delta}{\alpha-2}<0$, we choose a appropriate $\delta$ to determine the ranges of $\alpha$.
 
 \textbf{Case 1.} $\frac{\alpha+\delta}{2+\delta}>1$, $\alpha\neq-\delta\neq 2$ $\Longrightarrow$ $\alpha>2$ and $-2<\delta<0$. Then, from $Q<4\frac{\alpha+\delta}{\alpha-2}$, we have
 \[\frac{Q-4}{4}\left(\alpha-\frac{2Q}{Q-4}\right)<\delta<0,\]
 and
 \[\alpha<\frac{2Q}{Q-4}.\]
 Note that $\frac{2Q}{Q-4}>2$ for $Q>4$. In this case, we choose $\delta$ to satisfy $\max\{-2,\frac{Q-4}{4}(\alpha-\frac{2Q}{Q-4})\}<\delta<0$, which leads to $2<\alpha<\frac{2Q}{Q-4}$.
 
 \textbf{Case 2.} $\frac{\alpha+\delta}{2+\delta}>1$, $\alpha\neq-\delta\neq 2$ $\Longrightarrow$ $\alpha<2$ and $\delta<-2$. Under the assumption $Q<4\frac{\alpha+\delta}{\alpha-2}$, it follows that
 \[\delta<\frac{Q-4}{4}\left(\alpha-\frac{2Q}{Q-4}\right).\]
 Therefore, we choose $\delta$ such that $\delta<\min\{-2,\frac{Q-4}{4}(\alpha-\frac{2Q}{Q-4})\}$, which ensure that $\alpha<2$.
 
 We now examine the remaining cases, $\alpha=2$ and $\alpha=\frac{2Q}{Q-4}$.
 
 \textbf{Case 3.} $\alpha=2$. Returning to \eqref{Inequality-3}, let $\delta=-\alpha=-2$. Applying Young's inequality with the exponent pair $(\frac{\theta}{\theta-4},\frac{\theta}{4})$ (choose $\theta>4$), we get
 \[2\int\eta^\theta\leq\varepsilon_7\int\eta^\theta+C(\varepsilon_7)R^{Q-\theta}.\]
 Take $\varepsilon_7$ small enough, $\theta>Q$ and let $R\rightarrow+\infty$, we derive the contradiction.
 
 \textbf{Case 4.} $\alpha=\frac{2Q}{Q-4}$. In this case, we multiply both side of equation \eqref{Equation-sigma_2} by $2\eta^\theta$ and integrate over $\mathbb{H}^n$. Following exactly the same procedure as in the derivation of \eqref{Inequality-1+2+3}, we obtain
 \begin{equation}\label{Inequality-5}
  \begin{aligned}
   2\int(-u)^\alpha\eta^\theta
   \leq&-24n\int(Tu)^2\eta^\theta+12|\theta|\frac{C(Q)}{R}\int|Tu||\partial u|\eta^{\theta-1}\\
   &+\frac{3}{4}|\theta|(|\theta-1|+1)\frac{C(Q)}{R^2}\int|\partial u|^2\eta^{\theta-2}\\
   &+\frac{3}{2}|\theta|\frac{C(Q)}{R}\int\sigma_1|\partial u|\eta^{\theta-1}.
  \end{aligned}
 \end{equation}
 In fact, this can be recovered from \eqref{Inequality-1+2+3} simply by setting $\delta=0$. Recalling the estimate for \circled{5}, we readily have
 \begin{align*}
  12|\theta|\frac{C(Q)}{R}\int|Tu||\partial u|\eta^{\theta-1}\leq12|\theta|C(Q)\varepsilon_2\int(Tu)^2\eta^\theta+\frac{3}{\varepsilon_2}|\theta|\frac{C(Q)}{R^2}\int|\partial u|^2\eta^{\theta-2}.
 \end{align*}
 Substituting this into \eqref{Inequality-5} and choosing $\varepsilon_2$ sufficiently small, we arrive at
 \begin{equation}\label{Inequality-6}
  2\int(-u)^\alpha\eta^\theta\leq C_1R^{-1}\int\sigma_1|\partial u|\eta^{\theta-1}+C_2R^{-2}\int|\partial u|^2\eta^{\theta-2}\stackrel{\triangle}{=}\circled{8}+\circled{9},
 \end{equation} 
 where $C_1$ and $C_2$ are positive constants depending only on $Q$ and $\theta$, with $\theta$ to be determined.
 
 Let us now estimate the terms \circled{8} and \circled{9} more carefully. For the term \circled{8}, we begin by applying Young's inequality with the exponent pair $(2,2)$, then perform integration by parts, and finally use Young's inequality twice again, with the exponent pairs $(2,2)$ and $(4,\frac{4}{3})$. The detailed computation is given by
 \begin{align*}
  \circled{8}
  \leq&~\varepsilon_8\int\sigma_1|\partial u|^2(-u)^{\gamma-1}\eta^\theta+\frac{C_1^2}{4\varepsilon_8}R^{-2}\int\sigma_1(-u)^{-\gamma+1}\eta^{\theta-2}\\
  \leq&~\varepsilon_8\int\sigma_1|\partial u|^2(-u)^{\gamma-1}\eta^\theta\\
  &+\frac{C_1^2}{4\varepsilon_8}|1-\gamma|R^{-2}\int|\partial u|^2(-u)^{-\gamma}\eta^{\theta-2}+\frac{C_1^2}{4\varepsilon_8}|2-\theta|R^{-3}\int|\partial u|(-u)^{-\gamma+1}\eta^{\theta-3}\\
  \leq&~\varepsilon_8\int\sigma_1|\partial u|^2(-u)^{\gamma-1}\eta^\theta+(\varepsilon_9+\varepsilon_{10})\int|\partial u|^4(-u)^{\gamma-2}\eta^\theta\\
  &+\frac{C_1^4}{64}|1-\gamma|^2\varepsilon_8^{-2}\varepsilon_9^{-1}R^{-4}\int(-u)^{-3\gamma+2}\eta^{\theta-4}\\
  &+3\left(\frac{C_1}{4}\right)^{\frac{8}{3}}|2-\theta|^{\frac{4}{3}}\varepsilon_8^{-\frac{4}{3}}\varepsilon_{10}^{-\frac{1}{3}}R^{-4}\int(-u)^{-\frac{5}{3}\gamma+2}\eta^{\theta-4}.
 \end{align*}
 Here $\gamma<0$ is a parameter to be choosen later. The term \circled{9} can be estimated in the same way as \circled{8}, that is,
 \[\circled{9}\leq\varepsilon_{11}\int|\partial u|^4(-u)^{\gamma-2}\eta^\theta+\frac{C_2^2}{4\varepsilon_{11}}R^{-4}\int(-u)^{-\gamma+2}\eta^{\theta-4}.\]
 Plugging the above estimates into \eqref{Inequality-6} yields
 \begin{equation}\label{Inequality-7}
  \begin{aligned}
   2\int(-u)^\alpha\eta^\theta
   \leq&~\varepsilon_8\int\sigma_1|\partial u|^2(-u)^{\gamma-1}\eta^\theta+(\varepsilon_9+\varepsilon_{10}+\varepsilon_{11})\int|\partial u|^4(-u)^{\gamma-2}\eta^\theta\\
   &+\frac{C_1^4}{64}|1-\gamma|^2\varepsilon_8^{-2}\varepsilon_9^{-1}R^{-4}\int(-u)^{-3\gamma+2}\eta^{\theta-4}\\
   &+3\left(\frac{C_1}{4}\right)^{\frac{8}{3}}|2-\theta|^{\frac{4}{3}}\varepsilon_8^{-\frac{4}{3}}\varepsilon_{10}^{-\frac{1}{3}}R^{-4}\int(-u)^{-\frac{5}{3}\gamma+2}\eta^{\theta-4}\\
   &+\frac{C_2^2}{4}\varepsilon_{11}^{-1}R^{-4}\int(-u)^{-\gamma+2}\eta^{\theta-4}.
  \end{aligned}
 \end{equation}
 
 Setting $\alpha=\gamma<0$ in view of \eqref{Inequality-2}, and choosing $\varepsilon_i$ ($i=1,2,\cdots,5$) sufficiently small, we obtain
 \begin{equation}\label{Inequality-8}
  \int\sigma_1|\partial u|^2(-u)^{\gamma-1}\eta^\theta+\int|\partial u|^4(-u)^{\gamma-2}\eta^\theta\leq C_3R^{-4}\int(-u)^{\gamma+2}\eta^{\theta-4}.
 \end{equation}
 Therefore, making use of \eqref{Inequality-8} in \eqref{Inequality-7}, we derive the integral estimate
 \begin{equation}\label{Inequality-9}
  \begin{aligned}
   2\int(-u)^\alpha\eta^\theta
   \leq&~C_3(\varepsilon_8+\varepsilon_9+\varepsilon_{10}+\varepsilon_{11})R^{-4}\int(-u)^{\gamma+2}\eta^{\theta-4}\\
   &+C_4\varepsilon_8^{-2}\varepsilon_9^{-1}R^{-4}\int(-u)^{-3\gamma+2}\eta^{\theta-4}\\
   &+C_5\varepsilon_8^{-\frac{4}{3}}\varepsilon_{10}^{-\frac{1}{3}}R^{-4}\int(-u)^{-\frac{5}{3}\gamma+2}\eta^{\theta-4}\\
   &+C_6\varepsilon_{11}^{-1}R^{-4}\int(-u)^{-\gamma+2}\eta^{\theta-4},
  \end{aligned}
 \end{equation}
 where $C_4=\frac{C_1^4}{64}|1-\gamma|^2$, $C_5=3\left(\frac{C_1}{4}\right)^{\frac{8}{3}}|2-\theta|^{\frac{4}{3}}$, and $C_6=\frac{C_2^2}{4}$. 
 
 In order to obtain the desired estimate (i.e., to apply H\"{o}lder's inequality), we require $-2<\gamma<0$ to satisfy
 \[\begin{cases}
  \gamma+2<\alpha.\\
  -3\gamma+2<\alpha,\\
  -\frac{5}{3}\gamma+2<\alpha,\\
  -\gamma+2<\alpha,
 \end{cases}\]
 which means the range $\max\{-2,-\frac{\alpha-2}{3}\}<\gamma<0$. In this case, where $\alpha=\frac{2Q}{Q-4}$, this leads to the condition $-\frac{8}{3(Q-4)}<\gamma<0$. We further take $\varepsilon_8=\varepsilon_9=\varepsilon_{10}=\varepsilon_{11}=R^{-s}$, where $s>0$ will be determined later. We now apply H\"{o}lder's inequality with the exponent pairs $(\frac{\alpha}{\gamma+2},\frac{\alpha}{\alpha-\gamma-2})$, $(\frac{\alpha}{2-3\gamma},\frac{\alpha}{\alpha-2+3\gamma})$, $(\frac{\alpha}{2-\frac{5}{3}\gamma},\frac{\alpha}{\alpha-2+\frac{5}{3}\gamma})$, and $(\frac{\alpha}{2-\gamma},\frac{\alpha}{\alpha-2+\gamma})$ to the four terms on the right hand side of \eqref{Inequality-9}, respectively. By choosing $\theta$ sufficiently large, we then derive the integral estimate
 \begin{align*}
  2\int(-u)^\alpha\eta^\theta
  \leq&~C_3R^{-s-4+Q\frac{\alpha-\gamma-2}{\alpha}}\left(\int(-u)^\alpha\eta^\theta\right)^{\frac{\gamma+2}{\alpha}}+C_4R^{3s-4+Q\frac{\alpha-2+3\gamma}{\alpha}}\left(\int(-u)^\alpha\eta^\theta\right)^{\frac{2-3\gamma}{\alpha}}\\
  &+C_5R^{\frac{5}{3}s-4+Q\frac{\alpha-2+\frac{5}{3}\gamma}{\alpha}}\left(\int(-u)^\alpha\eta^\theta\right)^{\frac{2-\frac{5}{3}\gamma}{\alpha}}+C_6R^{s-4+Q\frac{\alpha-2+\gamma}{\alpha}}\left(\int(-u)^\alpha\eta^\theta\right)^{\frac{2-\gamma}{\alpha}}.
 \end{align*}
 Referring back to \eqref{Inequality-6}, all terms on the right hand side are integrated over the region $\mathrm{supp}(\partial\eta)=B_{2R}\setminus B_R$. Once we choose $s=-\frac{\gamma Q}{\alpha}=-\frac{\gamma(Q-4)}{2}$, this leads to
 \begin{equation}\label{Inequality-10}
  \begin{aligned}
   \int_{\mathbb{H}^n}(-u)^\alpha\eta^\theta
   \leq&~C_3\left(\int_{B_{2R}\setminus B_R}(-u)^\alpha\eta^\theta\right)^{\frac{\gamma+2}{\alpha}}+C_4\left(\int_{B_{2R}\setminus B_R}(-u)^\alpha\eta^\theta\right)^{\frac{2-3\gamma}{\alpha}}\\
   &+C_5\left(\int_{B_{2R}\setminus B_R}(-u)^\alpha\eta^\theta\right)^{\frac{2-\frac{5}{3}\gamma}{\alpha}}+C_6\left(\int_{B_{2R}\setminus B_R}(-u)^\alpha\eta^\theta\right)^{\frac{2-\gamma}{\alpha}}.
  \end{aligned}
 \end{equation}
 Since $\frac{\gamma+2}{\alpha}$, $\frac{2-3\gamma}{\alpha}$, $\frac{2-\frac{5}{3}\gamma}{\alpha}$ and $\frac{2-\gamma}{\alpha}\in(0,1)$, it follows from \eqref{Inequality-10} that the integral $\int_{\mathbb{H}^n}(-u)^\alpha\eta^\theta$ is bounded. Therefore, we obtain that
 \[\int_{B_{2R}\setminus R_R}(-u)^\alpha\eta^\theta\longrightarrow 0\]
 as $R\rightarrow+\infty$. Returning to \eqref{Inequality-10} and letting $R\rightarrow+\infty$, we conclude that
 \[\int_{\mathbb{H}^n}(-u)^\alpha\eta^\theta\leq 0,\]
 which is clearly in contradiction with the assumption that $u<0$.
 
 Thus, the proof of Theorem \ref{Main Theorem} is now complete.
\end{proof}

\bibliographystyle{plain}

\bibliography{mybibliography-2-Hessian}

\begin{thebibliography}{10}

\bibitem{Birindelli-Capuzzo-Cutri-1997}
I.~Birindelli, I.~Capuzzo~Dolcetta, and A.~Cutr\`i.
\newblock Liouville theorems for semilinear equations on the {H}eisenberg
  group.
\newblock {\em Ann. Inst. H. Poincar\'e{} C Anal. Non Lin\'eaire},
  14(3):295--308, 1997.

\bibitem{Caffarelli-Nirenberg-Spruck-1985}
L.~Caffarelli, L.~Nirenberg, and J.~Spruck.
\newblock The {D}irichlet problem for nonlinear second-order elliptic
  equations. {III}. {F}unctions of the eigenvalues of the {H}essian.
\newblock {\em Acta Math.}, 155(3-4):261--301, 1985.

\bibitem{Caffarelli-Gidas-Spruck-1989}
Luis~A. Caffarelli, Basilis Gidas, and Joel Spruck.
\newblock Asymptotic symmetry and local behavior of semilinear elliptic
  equations with critical {S}obolev growth.
\newblock {\em Comm. Pure Appl. Math.}, 42(3):271--297, 1989.

\bibitem{Catino-Li-Monticelli-Roncoroni-2023}
Giovanni Catino, Yanyan Li, Dario~D Monticelli, and Alberto Roncoroni.
\newblock A {L}iouville theorem in the {H}eisenberg group.
\newblock {\em arXiv preprint arXiv:2310.10469}, 2023.

\bibitem{Catino-Monticelli-Roncoroni-2023}
Giovanni Catino, Dario~D. Monticelli, and Alberto Roncoroni.
\newblock On the critical {$p$}-{L}aplace equation.
\newblock {\em Adv. Math.}, 433:Paper No. 109331, 38, 2023.

\bibitem{Chang-Gursky-Yang-2003}
Sun-Yung~A. Chang, Matthew~J. Gursky, and Paul~C. Yang.
\newblock Entire solutions of a fully nonlinear equation.
\newblock In {\em Lectures on partial differential equations}, volume~2 of {\em
  New Stud. Adv. Math.}, pages 43--60. Int. Press, Somerville, MA, 2003.

\bibitem{Chen-Li-1991}
Wen~Xiong Chen and Congming Li.
\newblock Classification of solutions of some nonlinear elliptic equations.
\newblock {\em Duke Math. J.}, 63(3):615--622, 1991.

\bibitem{Ciraolo-Figalli-Roncoroni-2020}
Giulio Ciraolo, Alessio Figalli, and Alberto Roncoroni.
\newblock Symmetry results for critical anisotropic {$p$}-{L}aplacian equations
  in convex cones.
\newblock {\em Geom. Funct. Anal.}, 30(3):770--803, 2020.

\bibitem{Danielli-Garofalo-Nhieu-2003}
Donatella Danielli, Nicola Garofalo, and Duy-Minh Nhieu.
\newblock Notions of convexity in {C}arnot groups.
\newblock {\em Comm. Anal. Geom.}, 11(2):263--341, 2003.

\bibitem{Flynn-Vetois-2023}
Joshua Flynn and J{\'e}r{\^o}me V{\'e}tois.
\newblock Liouville-type results for the {CR} {Y}amabe equation in the
  {H}eisenberg group.
\newblock {\em arXiv preprint arXiv:2310.14048}, 2023.

\bibitem{Gao-Shi-Zhang-2025}
Zhenghuan Gao, Shujun Shi, and Yuzhou Zhang.
\newblock On the solvability for a pk-hessian inequality.
\newblock {\em arXiv preprint arXiv:2503.00895}, 2025.

\bibitem{Gidas-Spruck-1981}
B.~Gidas and J.~Spruck.
\newblock Global and local behavior of positive solutions of nonlinear elliptic
  equations.
\newblock {\em Comm. Pure Appl. Math.}, 34(4):525--598, 1981.

\bibitem{Jerison-Lee-1987}
David Jerison and John~M. Lee.
\newblock The {Y}amabe problem on {CR} manifolds.
\newblock {\em J. Differential Geom.}, 25(2):167--197, 1987.

\bibitem{Jerison-Lee-1988}
David Jerison and John~M. Lee.
\newblock Extremals for the {S}obolev inequality on the {H}eisenberg group and
  the {CR} {Y}amabe problem.
\newblock {\em J. Amer. Math. Soc.}, 1(1):1--13, 1988.

\bibitem{Jerison-Lee-1989}
David Jerison and John~M. Lee.
\newblock Intrinsic {CR} normal coordinates and the {CR} {Y}amabe problem.
\newblock {\em J. Differential Geom.}, 29(2):303--343, 1989.

\bibitem{Lu-Manfredi-Stroffolini-2004}
Guozhen Lu, Juan~J. Manfredi, and Bianca Stroffolini.
\newblock Convex functions on the {H}eisenberg group.
\newblock {\em Calc. Var. Partial Differential Equations}, 19(1):1--22, 2004.

\bibitem{Ma-Ou-2023}
Xi-Nan Ma and Qianzhong Ou.
\newblock A {L}iouville theorem for a class semilinear elliptic equations on
  the {H}eisenberg group.
\newblock {\em Adv. Math.}, 413:Paper No. 108851, 20, 2023.

\bibitem{Obata-1971}
Morio Obata.
\newblock The conjectures on conformal transformations of {R}iemannian
  manifolds.
\newblock {\em J. Differential Geometry}, 6:247--258, 1971/72.

\bibitem{Ou-2010}
Qianzhong Ou.
\newblock Nonexistence results for {H}essian inequality.
\newblock {\em Methods Appl. Anal.}, 17(2):213--223, 2010.

\bibitem{Ou-2013}
Qianzhong Ou.
\newblock Singularities and {L}iouville theorems for some special conformal
  {H}essian equations.
\newblock {\em Pacific J. Math.}, 266(1):117--128, 2013.

\bibitem{Ou-2025}
Qianzhong Ou.
\newblock On the classification of entire solutions to the critical p-{L}aplace
  equation.
\newblock {\em Math. Ann.}, 392(2):1711--1729, 2025.

\bibitem{Phuc-Verbitsky-2006}
Nguyen~Cong Phuc and Igor~E. Verbitsky.
\newblock Local integral estimates and removable singularities for quasilinear
  and {H}essian equations with nonlinear source terms.
\newblock {\em Comm. Partial Differential Equations}, 31(10-12):1779--1791,
  2006.

\bibitem{Phuc-Verbitsky-2008}
Nguyen~Cong Phuc and Igor~E. Verbitsky.
\newblock Quasilinear and {H}essian equations of {L}ane-{E}mden type.
\newblock {\em Ann. of Math. (2)}, 168(3):859--914, 2008.

\bibitem{Sciunzi-2016}
Berardino Sciunzi.
\newblock Classification of positive {$D^{1,p}(\mathbb{R}^N)$}-solutions to the
  critical {$p$}-{L}aplace equation in {$\Bbb{R}^N$}.
\newblock {\em Adv. Math.}, 291:12--23, 2016.

\bibitem{Serrin-Zou-2002}
James Serrin and Henghui Zou.
\newblock Cauchy-{L}iouville and universal boundedness theorems for quasilinear
  elliptic equations and inequalities.
\newblock {\em Acta Math.}, 189(1):79--142, 2002.

\bibitem{Stein-1993}
Elias~M. Stein.
\newblock {\em Harmonic analysis: real-variable methods, orthogonality, and
  oscillatory integrals}, volume~43 of {\em Princeton Mathematical Series}.
\newblock Princeton University Press, Princeton, NJ, 1993.
\newblock With the assistance of Timothy S. Murphy, Monographs in Harmonic
  Analysis, III.

\bibitem{Trudinger-Wang-1997}
Neil~S. Trudinger and Xu-Jia Wang.
\newblock Hessian measures. {I}.
\newblock volume~10, pages 225--239. 1997.
\newblock Dedicated to Olga Ladyzhenskaya.

\bibitem{Trudinger-Wang-1999}
Neil~S. Trudinger and Xu-Jia Wang.
\newblock Hessian measures. {II}.
\newblock {\em Ann. of Math. (2)}, 150(2):579--604, 1999.

\bibitem{Trudinger-Wang-2002}
Neil~S. Trudinger and Xu-Jia Wang.
\newblock Hessian measures. {III}.
\newblock {\em J. Funct. Anal.}, 193(1):1--23, 2002.

\bibitem{Trudinger-Zhang-2013}
Neil~S. Trudinger and Wei Zhang.
\newblock Hessian measures on the {H}eisenberg group.
\newblock {\em J. Funct. Anal.}, 264(10):2335--2355, 2013.

\bibitem{Vetois-2016}
J\'er\^ome V\'etois.
\newblock A priori estimates and application to the symmetry of solutions for
  critical {$p$}-{L}aplace equations.
\newblock {\em J. Differential Equations}, 260(1):149--161, 2016.

\bibitem{Wei-Wu-2024}
Juncheng Wei and Ke~Wu.
\newblock A {L}iouville theorem for superlinear parabolic equations on the
  {H}eisenberg group.
\newblock {\em Adv. Nonlinear Stud.}, 24(1):189--205, 2024.

\bibitem{Xu-2009}
Lu~Xu.
\newblock Semi-linear {L}iouville theorems in the {H}eisenberg group via vector
  field methods.
\newblock {\em J. Differential Equations}, 247(10):2799--2820, 2009.

\end{thebibliography}

\end{document}